\tolerance=10000
\raggedbottom

\baselineskip=15pt
\parskip=1\jot

\def\sk{\vskip 3\jot}

\def\heading#1{\vskip3\jot{\noindent\bf #1}}
\def\label#1{{\noindent\it #1}}


\def\ref#1;#2;#3;#4;#5.{\item{[#1]} #2,#3,{\it #4},#5.}
\def\refinbook#1;#2;#3;#4;#5;#6.{\item{[#1]} #2, #3, #4, {\it #5},#6.} 
\def\refbook#1;#2;#3;#4.{\item{[#1]} #2,{\it #3},#4.}


\def\({\bigl(}
\def\){\bigr)}

\def\ga{\gamma}
\def\ze{\zeta}
\def\la{\lambda}
\def\si{\sigma}
\def\ps{\psi}

\def\Ga{\Gamma}
\def\De{\Delta}

\def\calS{{\cal S}}
\def\calW{{\cal W}}

\input xy
\xyoption{all}

\def\edge#1{\ar@{-}[#1]}

\def\text#1{{\rm #1}}
\def\state#1{*+++[o][F]{#1}}
\def\fstate#1{*++++[o][F=]{#1}}

\def\Ex{{\rm Ex}}
\def\Var{{\rm Var}}

\def\Stwo#1#2{\left\{#1\atop#2\right\}}
\def\Eul#1#2{\left\langle{#1\atop#2}\right\rangle}

\def\pwrt#1{{\partial\over\partial #1}}

\def\Li{{\rm Li}}

{
\pageno=0
\nopagenumbers
\rightline{\tt egp.pub.rerev.arxiv.tex}
\vskip1in

\centerline{\bf Stochastic Service Systems, Random Interval Graphs and Search Algorithms}
\vskip0.5in

\centerline{Patrick Eschenfeldt}
\centerline{\tt peschenfeldt@hmc.edu}
\sk

\centerline{Ben Gross}
\centerline{\tt bgross@hmc.edu}
\sk

\centerline{Nicholas Pippenger}
\centerline{\tt njp@math.hmc.edu}
\sk

\centerline{Department of Mathematics}
\centerline{Harvey Mudd College}
\centerline{1250 Dartmouth Avenue}
\centerline{Claremont, CA 91711}
\vskip0.5in

\noindent{\bf Abstract:}
We consider several stochastic service systems, and study the asymptotic behavior of the moments of various quantities that have application to models for random interval graphs and algorithms for searching for an idle server or empty waiting station.
In two cases the moments turn out to involve Lambert series for the generating functions for the sums of powers of divisors of positive integers.
For these cases we are able to obtain complete asymptotic expansions for the moments of the quantities in question.
\vskip0.5in
\leftline{{\bf Keywords:} Queueing theory, interval graphs, Lambert series, asymptotic expansions.}
\sk
\leftline{{\bf Subject Classification:} 60K26, 90B22}

\vfill\eject
}

\heading{1. Introduction}

Our goal in this paper is the study of stochastic service systems, with an eye to two applications: models for random interval graphs and the analysis of search algorithms.
The systems we study are traditionally designated $M/M/1$ (independent exponentially distributed interarrival times, independent exponentially distributed service times and a single server) and
$M/M/\infty$ (independent exponentially distributed interarrival times, independent exponentially distributed service times and infinitely many servers).

In a previous paper, Pippenger [P] introduced models for random interval graphs based on stochastic service systems, and analyzed, among others, a model based on the $M/M/\infty$ system.
For this system, {\it idle periods\/} (intervals of time during which all servers are idle) alternate with {\it busy periods\/} (intervals of time during which at least one server is busy).
A random graph can be constructed by considering a busy period, letting the vertices correspond to customers served during this busy period, and putting an edge between two vertices if the service intervals of the corresponding customers overlap.
Since edges correspond to intersections between intervals in a totally ordered set, the resulting graph is an interval graph, so this procedure yields a model for random interval graphs.
(Other models for random interval graphs have been considered by Scheinerman [S1, S2] and by
Godehardt and Jaworski [G2].)

Suppose that in the system $M/M/\infty$ customers arrive at rate $\la$ and are served at rate $1$.
Pippenger [P] showed that the number $N$ of vertices in the corresponding random graph (which corresponds to the number of customers served during the busy period) is such that the distribution of $N/e^\la$ tends to that of an exponential with mean $1$ as $\la\to\infty$.
Furthermore, the chromatic number $K$ of the graph (which for an interval graph equals the number of vertices in the largest clique of the graph, and corresponds to the largest number of customers simultaneously in the system during the busy period) is such that $K/e\la$ tends to $1$ in probability as $\la\to\infty$.

Our first goal in this paper is to study the corresponding random interval graph model for the $M/M/1$ system.
In this case we must have $\la<1$ to ensure that the busy period is finite with probability one, and that the number of customers served during the busy period has finite expectation,
and we shall be interested in asymptotics as $\la\to 1$.
When there is only one server, customers who arrive when the server is busy must wait for service,
and a {\it service discipline\/} (which determines which of the waiting customers will be served next when the current service interval ends) must be specified.
As a result, the corresponding interval graph will depend on the service discipline used.
Consider, for example, a busy period consisting of six service intervals, with two new customers arriving during the first interval and one new customer arriving during each of the second, third and fourth intervals.
Then if the service discipline is ``first-come-first-served'', the resulting graph will be
$$
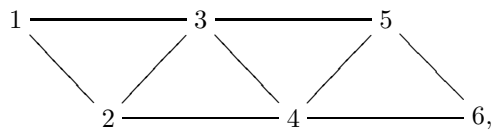
\xymatrix{
1\edge{rr} &&3\edge{rr} && 5 \\
& 2\edge{rr}\edge{ul}\edge{ur} &&4\edge{rr}\edge{ul}\edge{ur} && 6,\edge{ul} \\
}$$
whereas if the service discipline is ``last-come-first served'', the resulting graph will be
$$
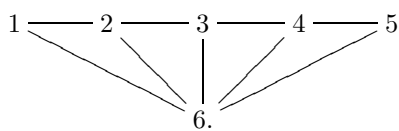
\xymatrix{
1\ar@{-}[r] &2\ar@{-}[r] &3\ar@{-}[r] &4\ar@{-}[r] &5 \\
&& 6.\ar@{-}[ull] \ar@{-}[ul] \ar@{-}[u] \ar@{-}[ur] \ar@{-}[urr] \\
}$$
Nevertheless, the number of vertices and the chromatic number the resulting graph 
(in this example, $6$ and $3$, respectively) are independent of the service discipline.
In Section 2, we shall derive asymptotic expansions for the moments of these quantities.
The leading terms of these expansions are
$$\Ex[N^m] \sim {2^{m-1} \, (2m-3)!! \, \la^{m-1}\over (1-\la)^{2m-1}} \eqno(1.1)$$
for $m\ge 1$, where $(2m-3)!! = (2m-3)\cdot(2m-5)\cdots 3\cdot 1$ and $(-1)!! = 1$,
$$\Ex[K] \sim \log {1\over 1-\la}, \eqno(1.2)$$
and, for $m\ge 2$,
$$\Ex[K^m] \sim {m!\,\ze(m) \over (1-\la)^{m-1}}, \eqno(1.3)$$
where $\ze(m) = \sum_{n\ge 1} 1/n^m$ is the Riemann zeta function.
It will be noted that $\Ex[K]$ grows quite slowly as $\la\to 1$.
If the random variable $J$ denotes the number of customers in the system in equilibrium, then 
$\Ex[J] = \la/(1-\la)$, which grows much more rapidly
(see Cohen [C2], p.~181).
It may appear paradoxical that the maximum number of customers in the system grows more slowly than the mean number of customers, but it must be borne in mind that $\Ex[K]$ is an average over busy periods, whereas $\Ex[J]$ is an average over time.
Indeed, the majority of busy periods have $K=1$: after the arrival initiating the busy period, 
the next event determines whether $K=1$ (if that event is a service termination) or $K>1$
(if that event is another arrival).
Because  $\la<1$,  the former (with probability $1/(1-\la)$) is more likely than the latter (with probability $\la/(1-\la)$).
We also note that, since $\ze(2) = \pi^2/6$, we have 
$\Var[K] = \Ex[K^2] - \Ex[K]^2 \sim \pi^2/3(1-\la)$, which grows much more rapidly 
than $\Ex[K]$ (or even $\Ex[K]^2$).
For the chromatic number, corresponding to the maximum number of customers in the system simultaneously during the busy period, the results involve Lambert series that are generating functions for arithmetical functions arising in number theory: the sums of powers of the divisors of positive integers.
These are
$$S_l(\la) = \sum_{n\ge 1} \si_l(n) \, \la^n, \eqno(1.4)$$
where $\si_l(n)$ denotes the sum of the $l$-th powers of the divisors of $n$
(see Hardy and Wright [H, p.~239]).

Our second goal in this paper is to analyze search algorithms connected with stochastic service systems.
Consider an infinite sequence $\calS_1, \calS_2, \ldots$ of servers in an $M/M/\infty$ system.
Suppose that each newly arriving customer scans this sequence in order and engages the first currently idle server.
We are interested in the index $L$ of the server $\calS_L$ engaged by a newly arriving customer in equilibrium.
This system has been extensively studied by Newell [N2], who suggests that $L$ ``is approximately uniformly distributed over the interval'' $[1, \la]$,
basing this assertion on the approximation
$$\Pr[L>l] \approx 
\cases{
1 - \displaystyle{l\over\la}, &if $l<\la$, \cr
0, &if $l>\la$. \cr
} \eqno(1.5)$$
But no error bounds are given for this or other approximations stated by Newell, and not even 
the fact that the first moment has the asymptotic behavior
$$\Ex[L]\sim{{\la\over 2}} \eqno(1.6)$$
that it would have under the uniform distribution is established rigorously.
In Section 3 we shall give a rigorous version of (1.5) that will suffice to establish
not only (1.6), but also the next term,
$$\Ex[L] = {\la\over 2} + {1\over 2}\log \la + O(1), \eqno(1.7)$$
and more generally
$$\Ex[L^m] = {\la^m \over m+1} + {m\,\la^{m-1}\log\la\over 2} + O\left(\la^{m-1}\right) \eqno(1.8)$$
for $m\ge 1$.
In particular,
we have
$$\eqalign{
\Var[L] 
&= \Ex[L^2] - \Ex[L]^2 \cr
& = {\la^2\over 12} + {\la\log\la\over 2} + O(\la). \cr
}$$
Since the interval $[0,1]$ is bounded, formula (1.8) shows that the $m$-th moment of $L/\la$ tends to $1/(m+1)$ as $\la\to\infty$ for all $m\ge 1$, and thus suffices to show that the distribution of 
$L/\la$ tends to the uniform distribution on the interval  $[0,1]$.
We note that a problem that is in a sense dual to ours (finding the largest index of a busy server, rather than the smallest index of an idle server) has been treated by Coffman, Kadota and Shepp [C1].

Our final results concern the analogue of the preceding search problem for the $M/M/1$ system.
Here there is only a single server, but an infinite sequence $\calW_1, \calW_2, \ldots$ of 
{\it waiting stations}.
A customer arriving when the server is busy scans this sequence in order and waits at the first vacant station.
When the server becomes free and there is at least one customer waiting, it too scans this sequence in order, and serves the customer waiting at the first occupied station.
We are interested in the index $I$ of the station $\calW_I$ at which a newly arriving customer waits in equilibrium (taken to be zero if the server is idle at the time of arrival).
(The index of the first station which the server finds  occupied (taken to be zero if the service interval initiates a busy period) has, of course, the same distribution as $I$.)
We shall show that the distribution of the random variable $I$ is closely related to that of the random variable $K$ studied in Section 2, with $\Ex[I^m] = \la\,\Ex[K^m]$.
This fact allows asymptotic expansions for the moments of $I$ to be obtained from those of the moments of $K$ in a straightforward way, with the result that the leading terms are the same.
\sk

\heading{2. Random Interval Graphs}

Our goal in this section is to determine asymptotic expansions for
the moments of the size (number of vertices) and chromatic number (number of vertices in the largest clique) for the random interval graph corresponding to the $M/M/1$ system.
These quantities correspond to the number $N$ of customers served during the busy period and the maximum number $K$ of customers in the system simultaneously during the busy period.

The random variable $N$ has a Catalan distribution:
$$Pr[N=n] = {1\over 2n-1} \, {2n-1\choose n} p^{n-1} \, q^n, \eqno(2.1)$$
where $p = \la/(1+\la)$ and $q = 1/(1+\la)$ (see for example Cohen [C2, 
pp.~190--191], or Riordan [R1, pp.~64--65]).
This distribution can be derived as follows.
Let $J$ denote the number of customers in the system.
When the busy period begins, $J=1$.
During the busy period, $J$ is incremented whenever a new customer arrives, and $J$ is decremented whenever a service interval ends and a customer departs, until $J=0$, at which time the busy period ends (see Figure 1).
\medskip
$$\xymatrix{
&&\ar[d]^-{\text{start}} \\
 & \fstate{\text{\scriptstyle{stop}}}  & \state{1} \ar@/^/[r]^-p \ar[l]^q&\cdots \ar@/^/[r]^-p \ar@/^/[l]^-q & \state{j} \ar@/^/[r]^-p \ar@/^/[l]^-q & \cdots \ar@/^/[l]^-q
}$$
\medskip
\centerline{
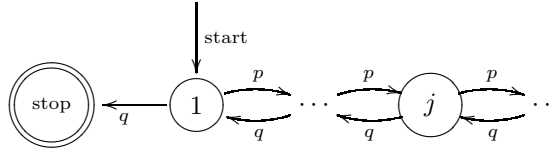
Figure 1. State transition diagram for determining} 
\centerline{the number of customers served during the busy period in 
$M/M/1$.
}
\medskip
When $J\ge 1$, the probability that the next transition is an arrival is $p = \la/(1+\la)$, and the probability that the next transition is a departure is $q = 1/(1+\la)$.
If $n$ customers are served during the busy period, there must be $n-1$ further arrivals 
(beyond the one that began the busy period) and
$n$ departures, and these must occur in such an order that $J=0$ for the first time immediately after the last of these $n$ departures.
The number of such orders is $A_n = {2n-1\choose n}/(2n-1) = {2n-2\choose n-1}/n$,
the $n$-th Catalan number (see for example Comtet [C3, p.~53]).
Thus we have (2.1).
Since the generating function for the Catalan numbers is
$a(z) = \sum_{n\ge 1} A_n\,z^n = \(1 - \sqrt{1-4z}\)/2$, the probability generating function $g_N(z)$ for $N$ is 
$$\eqalign{
g_N(z)
&= \sum_{n\ge 1} \Pr[N=n] \, z^n \cr
&= {a(pqz)\over p} \cr
&= {1 - \sqrt{1 - 4pqz} \over 2p}. \cr
}$$
Since for $m\ge 1$, we have $d^m a(z)/dz^m = 2^{m-1} \, (2m-3)!! / (1-4z)^{(2m-1)/2}$,
the factorial moments of $N$ are given by
$$\eqalignno{
\Ex[N(N-1)\cdots(N-m+1)]
&= {1\over p} {d^m \over dz^m} {1 - \sqrt{1-4pqz} \over 2}\bigg\vert_{z=1} \cr
&= {p^{m-1} \, q^m \, 2^{m-1} \, (2m-3)!! \over \sqrt{1 - 4pq}^{(2m-1)}} \cr
&= {2^{m-1} \, (2m-3)!! \, \la^{m-1} \over (1-\la)^{2m-1}}, \cr
}$$
because $\sqrt{1-4pq} = q-p = (1-\la)/(1+\la)$.
Since $x^m = \sum_{0\le l\le m} \Stwo{m}{l} \, x(x-1)\cdots(x-l+1)$, where the $\Stwo{m}{l}$
are the Stirling numbers of the second kind, with the generating function
$\sum_{m\ge l\ge 0} \Stwo{m}{l} {y^l \, z^m \over l!} = e^{y(e^z-1)}$
(see for example Comtet [C3, pp.~206--207]), and $\Stwo{m}{0} = 0$ for $m\ge 1$,
the ordinary moments of $N$ are given by
$$\eqalignno{
\Ex[N^m]
&= \sum_{1\le l\le m} \Stwo{m}{l} \, \Ex[N(N-1)\cdots(N-l+1)] \cr
&= \sum_{1\le l\le m} \Stwo{m}{l} \, {2^{l-1} \, (2l-3)!! \, \la^{l-1} \over (1-\la)^{2l-1}}. \cr
}$$
Since $\Stwo{m}{m} = 1$, this gives the asymptotic formula (1.1) as $\la\to 1$.

We turn now to the random variable $K$, the maximum number of customers in the system simultaneously during the busy period (counting the customer being served, so that $K\ge 1$).
The distribution of $K$ is given by
$$\Pr[K>k] = {(1-\la) \, \la^k \over 1 - \la^{k+1}} \eqno(2.2)$$
(see for example Cohen [C2, pp.~191--193]).
This distribution can be derived as follows.
Consider a game played between two players: $P$, who begins with $v$ dollars, and $Q$ who begins with $w$ dollars.
At each step of the game, a biassed coin is tossed; $P$ wins with probability $p$, in which case 
$Q$ pays $P$ one dollar, and $Q$ wins with the complementary probability $q=1-p$, in which case 
$P$ pays $Q$ one dollar.
The game continues until one of the players is ruined (that is, has no money left).
It is known that (1) with probability one, either $P$ or $Q$ is eventually ruined, and (2),
if $p\not=q$, then the probability that $Q$ is ruined is
$$\Pr[Q\hbox{\ ruined}] = {(q/p)^v - 1 \over (q/p)^{v+w} - 1} \eqno(2.3)$$
(see for example Feller [F, p.~345]).

Now consider a busy period of the $M/M/1$ queue.
The successive events of arrivals and terminations of service intervals during the busy period correspond to steps in the game described above. 
The wealth of player $P$ will correspond to the number $J$ of customers in the system, so $v=1$.
An arrival will correspond to a win by player $P$, so $p=\la/(1+\la)$, and
the termination of a service interval will correspond to a win by player $Q$, so $q=1/(1+\la)$.
Suppose that player $Q$ begins with $w=k$ dollars.
Then the event $K>k$ will correspond to $Q$ being ruined.
Substituting these values in (2.3) yields (2.2) (see Figure 2).
\medskip
$$\xymatrix{
&&\ar[d]^-{\text{start}} \\
 & *+++[o][F=]{\scriptstyle K\le k}  & \state{1} \ar@/^/[r]^-p \ar[l]^-q&\cdots \ar@/^/[r]^-p \ar@/^/[l]^-q & \state{j} \ar@/^/[r]^-p \ar@/^/[l]^-q & \cdots \ar@/^/[r]^-p \ar@/^/[l]^-q &\state{k}\ar[r]^-p  \ar@/^/[l]^-q&*+++[o][F=]{\scriptstyle K>k}
}$$
\medskip
\centerline{Figure 2. State transition diagram for determining whether $K\le k$ or $K>k$.}
\medskip
This correspondence also shows what happens for $\la\ge 1$.
For $\la=1$ (in which case the busy period is finite with probability one, but its expected length is infinite), we have take $p=q=1/2$, and have
$$\Pr[Q\hbox{\ ruined}] = {v\over v+w}.$$
This result yields
$$\Pr[K>k] = {1\over k+1},$$
so that 
$$\Ex[K] = \sum_{l\ge 0} \Pr[K>k] \eqno(2.4)$$ 
diverges logarithmically.
Of course, for $\la>1$ (in which case the busy period is infinite with positive probability),
(2.3) shows that (2.4) diverges linearly.
\vfill\eject

Our next goal is to determine the moments of $K$:
$$\Ex[K^m] = \sum_{k\ge 0} k^m \, \Pr[K=k].$$
Writing 
$$\eqalign{
\De_m(k) 
&= k^m - (k-1)^m \cr
&= \sum_{0\le l\le m-1} {m\choose l}(-1)^{m-1-l} \, k^l \cr
}$$
for the backward differences of the $m$-th powers, and setting
$$T_m(\la) = \sum_{n\ge 1} {n^m \la^n \over 1 - \la^n}, \eqno(2.5)$$
summation by parts yields
$$\eqalignno{
\Ex[K^m]
&= \sum_{k\ge 0} k^m \, \Pr[K=k] \cr
&= \sum_{k\ge 0} \De_m(k+1) \,\Pr[K>k] \cr
&= (1-\la) \sum_{k\ge 0} {\De_m(k+1) \,\la^k \over 1-\la^{k+1}} \cr
&= {1-\la\over \la} \sum_{j\ge 1} {\De_m(j) \,\la^j \over 1-\la^j} \cr
&= {1-\la\over \la} \sum_{j\ge 1} \sum_{0\le l\le m-1} {m\choose l} \, (-1)^{m-1-l} 
\, {j^l \,\la^j \over 1-\la^j} \cr
&= {1-\la\over \la} \sum_{0\le l\le m-1} {m\choose l} \, (-1)^{m-1-l} \, T_l(\la). &(2.6)\cr
}$$
Since $\Ex[K^m]$ is a linear combination of the $T_l(\la)$, it will suffice to determine the asymptotic behavior of the sums $T_m(\la)$.
The sums $T_l(\la)$ are in fact the Lambert series $S_l(\la)$ given by (1.4); 
we have
$$\eqalignno{
T_l(\la)
&=\sum_{j\ge 1} {j^l \, \la^j \over 1-\la^j} \cr
&= \sum_{j\ge 1} j^l \, \sum_{i\ge 1} \la^{ij} \cr
&= \sum_{i\ge 1} \sum_{j\ge 1} j^l \, \la^{ij} \cr
&= \sum_{n\ge 1} \,\la^n \, \sum_{d\mid n} d^l &(2.7)\cr
&= \sum_{n\ge 1} \si_l(n)\,\la^n, \cr
&= S_l(\la) \cr
}$$
where the inner sum in (2.7) is over integers $d$ dividing $n$.

We note that the sums $T_l(\la)$ can be expressed in terms of known (albeit exotic) functions of analysis.
We define the $q$-gamma function by
$$\Ga_q(x) = (1-q)^{1-x} \prod_{n\ge 0} {1-q^{n+1}\over 1-q^{n+x}}$$
(see for example Gasper and Rahman[G1, p.~16]). 
(This function gets its name from the fact that
$\lim_{q\to 1} \Ga_q(x) = \Ga(x)$, where $\Ga(x)$ is the Euler gamma function; see for example
Whittaker and Watson [W, pp.~235--264].)
If we define the $q$-digamma function $\ps_q(x)$ as the logarithmic derivative
$$\eqalign{
\ps_q(x)
&= \pwrt{x} \log \Ga_q(x) \cr
&= -\log(1-q) + \log q \, \sum_{n\ge 0} {q^{n+x} \over 1-q^{n+x}} \cr
}$$
of the $q$-gamma function, then we have
$$T_0(\la) = {\ps_\la(1) + \log(1-\la) \over \log\la}.$$
To go further, we define the $l$-th $q$-polygamma function $\ps_q^{(l)}$ as the $l$-th derivative
$$\ps_q^{(l)}(x) = \left(\pwrt{x}\right)^l \ps_q(x)$$
of the $q$-digamma function.
If we set $z=q^{n+x}$, then
$$\left(z\pwrt{z}\right) = \left({1\over \log q}\,\pwrt{x}\right).$$
Since
$$\sum_{i\ge 1} i^l z^i = \left(z\pwrt{z}\right)^l {z\over 1-z},$$
we have
$$\sum_{i\ge 1} i^l \, q^{i(n+x)} 
= {1\over \log^l q} \,  \left(\pwrt{x}\right)^l \, {q^{n+x} \over 1 - q^{n+x}}.$$
Summing over $n\ge 0$ yields
$$\eqalign{
\sum_{i\ge 1} {i^l \, q^{ix} \over 1 - q^{ix}}
&= \sum_{i\ge 1} i^l \,  \sum_{n\ge 0}  q^{i(n+x)} \cr
&= {1\over \log^l q} \, \left(\pwrt{x}\right)^l \, 
\sum_{n\ge 0} {q^{n+x} \over 1 - q^{n+x}} \cr
&= {1\over \log^l q} \, \left(\pwrt{x}\right)^l \, 
 {\ps_q(x) + \log(1-q) \over \log q}. \cr
}$$
Thus for $l\ge 1$ we have
$$T_l(\la) = {\ps_\la^{(l)}(1) \over \log^{l+1} \la}.$$

Our next goal is to derive the leading terms (1.2) and (1.3) of the asymptotic expansion of the moments of $K$.
To establish (1.2) and (1.3), we begin by deriving
$$T_0(\la)  \sim {1\over 1-\la}\log {1\over 1-\la} \eqno(2.8)$$
and, for $l\ge 1$,
$$T_l(\la)  \sim {l! \, \ze(l+1)\over (1-\la)^{l+1}}. \eqno(2.9)$$
Once these formulas are established, it will be clear that the sum in (2.6) is dominated by the term for which $l=m-1$, so that $\Ex[K^m] \sim m\,S_{m-1}(\la)$, and (1.2) and (1.3) follow from (2.8) and (2.9), respectively.

Our strategy for proving (2.8) and (2.9) will be to approximate the sums
$S_l(\la)$
by integrals
$$I_l(\la) = \int_1^\infty {x^l \, \la^x \, dx \over 1-\la^x },$$
then then to show that the difference $S_l(\la)-I_l(\la)$ is negligible in comparison with $I_l(\la)$.
It will be convenient to write $\la = e^{-h}$.
The limit $\la\to 1$ then corresponds to $h\to 0$.
We have
$$\eqalignno{
h
&= \log{1\over \la} \cr 
&= \log{1\over 1-(1-\la)} \cr
&\sim 1-\la. &(2.10) \cr
}$$

For $l=0$, we have
$$\eqalign{
I_0(\la) 
&= \int_1^\infty { \la^x \, dx \over 1-\la^x } \cr
&= \int_1^\infty  \sum_{l\ge 1}  \, \la^{lx} \,dx  \cr
&=  \sum_{l\ge 1}  \, \int_1^\infty  e^{-hlx} \,dx  \cr
&=  \sum_{l\ge 1}  \, {e^{-hl} \over hl }  \cr
&=  {1\over h} \sum_{l\ge 1}  \, {\la^{l} \over l }  \cr
&= {1 \over h }\;\log {1\over 1-\la}. \cr
}$$
Substituting (2.10) in this result yields 
$$I_0(\la) \sim  {1\over 1-\la}\log {1\over 1-\la}. \eqno(2.11)$$
We bound $\vert S_0(\la) - I_0(\la)\vert$ by the total variation of $f(x) = \la^x / (1-\la^x)$.
Since $f(x)$ decreases monotonically from $\la/(1-\la)$ to $0$ as $x$ increases from 
$1$ to $\infty$, we have $\vert S_0(\la) - I_0(\la)\vert \le \la/(1-\la) \sim 1/(1-\la)$.
Since this difference is negligible in comparison with (2.11), we obtain (2.8).

For $l\ge 1$, we have
$$\eqalignno{
I_l(\la)
&= \int_1^\infty {x^l \, \la^x \, dx \over 1-\la^x} \cr
&= \int_0^\infty {(y+1)^l \, \la^{y+1} \, dy \over 1-\la^{y+1}} \cr
&= \int_0^\infty \sum_{0\le i\le l} {l\choose i}{y^i \, \la^{y+1} \, dy \over 1-\la^{y+1}} \cr
&= \int_0^\infty \sum_{0\le i\le l} {l\choose i}y^i \sum_{j\ge 1} \la^{j(y+1)} \, dy  \cr
&= \int_0^\infty \sum_{0\le i\le l} {l\choose i}y^i \sum_{j\ge 1} e^{-hj(y+1)} \, dy  \cr
&=  \sum_{0\le i\le l} {l\choose i}\sum_{j\ge 1}\int_0^\infty y^i \, e^{-hj(y+1)} \, dy  \cr
&=  \sum_{0\le i\le l} {l\choose i}\sum_{j\ge 1} {i! \, e^{-hj} \over (hj)^{i+1}} \cr
&= \sum_{0\le i\le l} {l\choose i} {i! \over h^{i+1}} \sum_{j\ge 1} {\la^{j} \over j^{i+1}} \cr
&= \sum_{0\le i\le l} {l\choose i} {i! \over h^{i+1}} \Li_{i+1}(\la), &(2.12)\cr
}$$
where $\Li_k(\la) = \sum_{n\ge 1} \la^n/n^k$ is the $k$-th polylogarithm.
Since $\Li_1(\la) = \log \(1/(1-\la)\)$ and $\Li_k(\la)\to\ze(k)$ as $\la\to 1$ for $k\ge 2$,
the sum in (2.12) is dominated by the term for which $i=l$, and we have
$$\eqalignno{
I_l(\la) 
&\sim {l! \, \ze(l+1) \over h^{l+1}} \cr
&\sim {l! \, \ze(l+1) \over (1-\la)^{l+1}}&(2.13) \cr
}$$
We bound $\vert S_l(\la)-I_l(\la)\vert$ by the total variation of $f(x) = x^l \, \la^x / (1-\la^x)$ for
$0\le x < \infty$.
As $x$ increases, $f(x)$ increases monotonically from $0$ to a maximum, then decreases
monotonically to $0$.
Thus the total variation of $f(x)$ is twice the maximum.
This maximum is
$$\eqalign{
\max_{0\le x < \infty} f(x)
&= \max_{0\le x < \infty} {x^l \, e^{-hx} \over 1 - e^{-hx}} \cr
&= \max_{0\le x < \infty} {x^l  \over e^{hx} - 1} \cr
&= {1\over h^l} \, \max_{0\le y < \infty}  {y^l \over e^y - 1}. \cr
}$$
Furthermore, $y^l / (e^y - 1) \le l!$, because $e^y - 1 = \sum_{n\ge 1} y^n/n! \ge y^l/l!$.
Thus $\vert S_l(\la)-I_l(\la)\vert \le 2\max_{0\le x<\infty} f(x) \le 2\,l!/h^l \sim 2\,l!/(1-\la)^l$.
Since this difference is negligible in comparison with (2.13), we obtain (2.9).

We shall now show how asymptotic expansions, with error terms of the form
$O\((1-\la)^R\)$ for any $R$, can be derived for all of the moments $\Ex[K^m]$.
The essence of the argument is to use the Euler-Maclaurin formula to estimate the difference
between $S_l(\la)$ and $I_l(\la)$.
This is most conveniently done using a result of Zagier [Z].
Indeed, for $l\ge 1$, Zagier gives the expansion for $S_l(\la)$, in terms of the parameter
$h = -\log\la$ rather than $1-\la$.
All that remains for us to do is substitute an expansion for $h$ in terms of $1-\la$.
For $l=0$, the expansion for $S_0(\la)$ in terms of $h$ has been given by Egger (n\'{e} Endres) and Steiner [E1, E2], again using the result of Zagier.
We shall proceed differently, to obtain an expansion involving $-\log(1-\la)$ rather than $-\log h$.

\label{Proposition:}
(Zagier [Z, p.~318])
Let $f(x)$ be analytic at $x=0$,
with power series $f(x) = \sum_{r\ge 0} b_r \, x^r$ about $x=0$.
Suppose that $\int_0^\infty \vert f^{(r)}(x)\vert \, dx < \infty$ for all $r\ge 0$,
where $f^{(r)}(x)$ denotes the $r$-th derivative of $f(x)$.
Define $F = \int_0^\infty f(x) \, dx$.
Let $g(x) = \sum_{n\ge 1} f(nx)$.
Then $g(x)$ has the asymptotic expansion
$$g(x) \sim {F \over x} + \sum_{r\ge 0} {b_r \, B_{r+1} \, (-1)^r x^r \over (r+1) }, \eqno(2.14)$$
where $B_r$ is the $r$-th Bernoulli number, defined by $t/(e^t - 1) = \sum_{r\ge 0} B_r \, t^r/r!$.

This result is proved by using the Euler-Maclauren formula,
$$\eqalign{
\int_0^N f(y)\,dy
&= {f(0)\over 2} + \sum_{1\le n\le N-1} f(n) + {f(N)\over 2}
+ \sum_{1\le r\le R-1} {(-1)^r B_{r+1} \over (r+1)!}\(f^{(r)}(N) - f^{(r)}(0)\) \cr
&\qquad +(-1)^R \int_0^N f^{(R)}y {{B}_R(\{y\})\over R!}\,dy, \cr
}$$
where ${B}_r(y)$ is the  $r$-th Bernoulli polynomial, defined by 
$t e^{yt}/(e^t - 1) = \sum_{r\ge 0} B_r(y) \, t^r/r!$, and $\{y\} = y-\lfloor y\rfloor$ denotes the fractional part of $y$ .
(For the Euler-Maclauren formula, the Bernoulli numbers and the Bernoulli polynomials,
see for example Whittaker and Watson [W, pp.~125--128], where, however, the indexing of the numbers and polynomials is different.)
The condition $\int_0^\infty \vert f^{(r)}(y)\vert \, dy < \infty$ allows us to let $N\to\infty$, obtaining
$$\eqalign{
\int_0^\infty f(y)\,dy 
&=  \sum_{n\ge 1} f(n) 
- \sum_{0\le r\le R-1} {(-1)^r B_{r+1} \over (r+1)!} \, f^{(r)}(0) 
 +(-1)^R \int_0^\infty f^{(R)}(y) {{B}_R(\{y\})\over R!}\,dy. \cr
}$$
If we now write $f(xy)$ instead of $f(y)$, we obtain
$$\eqalign{
\int_0^\infty f(xy)\,dy 
&=  \sum_{n\ge 1} f(nx) 
- \sum_{0\le r\le R-1} {(-1)^r B_{r+1} \over (r+1)!} \, f^{(r)}(0) \, x^r
 +(-1)^R \, x^R \int_0^\infty f^{(R)}(xy) {{B}_R(\{y\})\over R!}\,dy. \cr
}$$
Changing the variable of integration from $y$ to $y/x$ then yields
$$\eqalign{
{1\over x} \int_0^\infty f(y)\,dy 
&=  \sum_{n\ge 1} f(nx) 
- \sum_{0\le r\le R-1} {(-1)^r B_{r+1} \over (r+1)!} \, f^{(r)}(0) \, x^r
 +(-1)^R \, x^{R-1} \int_0^\infty f^{(R)}(y) {{B}_R(\{y/x\})\over R!}\,dy. \cr
}$$
The integral on the left-hand side is $F$, the first sum on the right-hand side is $g(x)$,
$f^{(r)}(0) = r!\,b_r$, and the last term on the right-hand side is $O(x^{R-1})$.
Thus
$${F\over x} = g(x) - \sum_{0\le r\le R-1} {b_r \, B_{r+1} \, (-1)^r x^r \over (r+1) }
+ O(x^{R-1}),$$
which yields the expansion (2.14).

For $l\ge 1$, we define
$$f(x) = {x^l \over e^x -1}.$$
Then $f(x)$ is analytic at $x=0$ with the Taylor series
$$f(x) = \sum_{r\ge 0} {B_r \, x^{r+l-1} \over r!}$$
and the integral
$$\eqalign{
F
&= \int_0^\infty {x^l \, e^{-x} \, dx \over 1 - e^{-x}} \cr
&= l! \, \ze(l+1) \cr
}$$
(see for example Whittaker and Watson [W, p.~266]).
Furthermore, $f^{(r)}(x)$ is a rational function of $x$ and $e^x$, in which the degree of the numerator in $e^x$ is $r$, while the denominator is $(e^x-1)^{r+1}$.
Thus $f(x)$ satisfies the conditions of the proposition, and
we have the asymptotic expansion
$$g(x) \sim {l! \, \ze(l+1) \over x} + 
\sum_{r\ge 0} {(-1)^{r+l-1} \, B_r \, B_{r+l} \, x^{r+l-1} \over r! \, (r+l)}.$$
Recalling that $\la = e^{-h}$, so that $h=-\log\la$, we therefore have
$$\eqalignno{
T_l(\la)
&= \sum_{n\ge 1} {n^l \, e^{-hn} \over 1 - e^{-hn}} \cr
&= {1\over h^l} \sum_{n\ge 1} {(nh)^l  \over e^{nh} - 1} \cr
&= {1\over h^l} \sum_{n\ge 1} f(nh) \cr
&= {g(h) \over h^l}  \cr
&\sim {l! \, \ze(l+1) \over h^{l+1}} + 
\sum_{r\ge 0} {(-1)^{r+l-1} \, B_r \, B_{r+l} \, h^{r-1} \over r! \, (r+l)}. &(2.15)\cr
}$$
We note that, if $l$ is odd, then this expansion has only finitely many terms (because $B_r=0$ for odd $r\ge 3$).
To obtain an asymptotic expansion in terms of $1-\la$, we must substitute the expansion for $1/h$:
$$\eqalignno{
{1\over h}
&= {1\over -\log \la} \cr
&= {-1\over \log \(1-(1-\la)\)} \cr
&= {1\over 1-\la}\, {-(1-\la)\over \log \(1-(1-\la)\)} \cr
&= {1\over 1-\la}\, \sum_{r\ge 0} {(-1)^r \, C_r \, (1-\la)^r \over r!}, &(2.16)\cr
}$$
where $C_r$ is the $r$-th Bernoulli number of the second kind, defined by
$t/\log(1+t) = \sum_{r\ge 0} C_r \, t^r/r!$ (see for example Roman [R2, p.~116]).
(These numbers are also called the Cauchy numbers of the first kind, and are given by
$C_r = \int_0^1 x(x-1)\cdots(x-r+1)\,dx$; see for example Comtet [C2, pp.~293--294].)

For $l=0$, we must proceed differently, because 
$$f(x) = {1\over e^x - 1}$$
has a pole at $x=0$.
We define
$$\eqalign{
f^*(x) 
&= f(x)  - {e^{-x} \over x} \cr
&= {1\over e^x - 1} - {e^{-x} \over x}. \cr
}$$
Then $f^*(x)$ is analytic at $x=0$ with the Taylor series
$$f^*(x) = \sum_{r\ge 0} {\(B_{r+1} - (-1)^{r+1}\) \, x^r \over (r+1)!}$$
and the integral
$$\eqalign{
F^*
&= \int_0^\infty  \left({1\over e^x - 1} - {e^{-x} \over x}\right) \, dx\cr
& = \ga \cr
}$$
(see for example Whittaker and Watson [W, p.~246]).
Furthermore, $f^{*(R)}(x)$ is a rational function of $x$ and $e^x$, in which the degree of the numerator in $e^x$ is $R$, while the denominator is $\((e^x-1)\,x\)^{R+1}$.
Thus $f^*(x)$ satisfies the conditions of the proposition, and
we have the asymptotic expansion
$$g^*(x) \sim {\ga\over x} 
+ \sum_{r\ge 0} {(-1)^r \, B_{r+1} \(B_{r+1} - (-1)^{r+1}\) \, x^r \over (r+1)\,(r+1)!}.$$
We therefore have
$$\eqalignno{
T_0(\la)
&= \sum_{n\ge 1} {e^{-nh} \over 1 - e^{-nh}} \cr
&= \sum_{n\ge 1} {1 \over e^{nh} - 1} \cr
&= \sum_{n\ge 1} {e^{-nh} \over nh} + \sum_{n\ge 1} {1 \over e^{nh} - 1} - {e^{-nh} \over nh} \cr
&= {1\over h}\log {1\over 1-\la} + \sum_{n\ge 1} f^*(nh) \cr
&= {1\over h}\log {1\over 1-\la} +  g^*(h) \cr
&\sim {1\over h}\log {1\over 1-\la} + {\ga\over h} 
+ \sum_{r\ge 0} {(-1)^r \, B_{r+1} \(B_{r+1} - (-1)^{r+1}\) \, h^r \over (r+1)\,(r+1)!}. &(2.17)\cr
}$$

To obtain asymptotic expansions for the moments of $K$, we substitute (2.16) into (2.15) and (2.17),
then substitute the results into (2.6),
using the expansion
$$\eqalign{
{1-\la\over\la}
&= {1-\la\over 1-(1-\la)} \cr
&= \sum_{r\ge 1} (1-\la)^r. \cr
}$$
Retaining only terms that do not vanish as $\la\to 1$, we obtain
$$\Ex[K] = \log{1\over 1-\la} + \ga + O\left((1-\la)\log{1\over 1-\la}\right)$$
and
$$\Ex[K^2] = {\pi^2\over 3}{1\over 1-\la} + \log{1\over 1-\la} 
+ (\ga-1) + O\left((1-\la)\log{1\over 1-\la}\right)$$
for the first two moments.
Thus we have
$$\eqalign{
\Var[K]
&= \Ex[K^2] - \Ex[K]^2 \cr
&= {\pi^2\over 3}{1\over 1-\la} - \log^2 {1\over 1-\la} + (1-2\ga)\log{1\over 1-\la} - \ga^2
+O\left((1-\la)\log^2 {1\over 1-\la}\right). \cr
}$$
\sk

\heading{3. Search Algorithms}

In this section we shall analyze search algorithms for $M/M/\infty$ and $M/M/1$ systems.
We begin by studying the distribution of the random variable $L$, defined as the index
of the first idle server $\calS$ found by an arriving customer in the $M/M/\infty$ system.
Our goal is to prove (1.8), which gives the first two terms in the asymptotic expansions of the moments of $L$.
The key to our results is the probability $\Pr[L>l]$, which is simply the probability that the first $l$ servers $\calS_1,\ldots, \calS_l$ are all busy.
It is well known that this probability is given by the Erlang loss formula
$$\eqalign{
\Pr[L>l] 
&= {\la^l/l! \over \sum_{0\le k\le l} \la^k/k!} \cr
&= {1\over D_l}, \cr
}$$
where
$$D_l = \sum_{0\le k\le l} {l! \over (l-k)! \, \la^k} \eqno(3.1)$$
(see for example Newell [N, p.~3]).
The sum $D_l$ can be expressed as an integral,
$$D_l = \int_0^\infty \left(1 + {x\over\la}\right)^l \, e^{-x} \, dx$$
(see for example Newell [N, p.~7]), and most of Newell's analysis is based on such a representation. But we shall work directly with the expression of $D_l$ as the sum in (3.1).
We shall partition the values of $l$ into two ranges.
The first, which we shall call the ``body'' of the distribution, will be $0\le l\le l_0 = \la-s$,
where $s = \sqrt{\la}$.
The second, which we shall call the ``tail'', will be $l > l_0$.

We begin with the body.
We shall establish the estimate
$$\Pr[L>l] = (1 - l/\la) +  {1 \over \la(1 - l/\la)}  + O\left({1\over\la}\right)
+ O\left({1\over \la^2(1-l/\la)^3}\right) \eqno(3.2)$$
for $l\le l_0 = \la - s$, where $s = \sqrt{\la}$.
We begin by using the principle of inclusion-exclusion to derive bounds on the denominator
$D_l$.

We begin with a lower bound.
Since
$$\eqalign{
l(l-1)\cdots(l-k+1) 
&\ge l^k  - \left(\sum_{0\le j\le k-1} j\right) l^{k-1} \cr
&= l^k - {k\choose 2}l^{k-1}, \cr
}$$
we have
$$\eqalign{
D_l
&= \sum_{0\le k\le l} {l(l-1)\cdots(l-k+1)\over\la^k} \cr
&\ge \sum_{0\le k\le l} \left({l\over\la}\right)^k  
- {1\over\la}\sum_{0\le k\le l} {k\choose 2}\left({l\over\la}\right)^{k-1}. \cr
}$$
For the first sum we have
$$\sum_{0\le k\le l} \left({l\over\la}\right)^k   = {1 + O\((l/\la)^l\)\over 1 - l/\la}.$$
We note that the logarithm of $(l/\la)^l$ has a non-negative second derivative for $l\ge 1$.
Thus $(l/\la)^l$ assumes its maximum in the interval $0\le l\le l_0$ for $l=0$,  $l=1$ or
$l=l_0$.
Its values there are $0$, $1/\la$ and $(1-s/\la)^{\la-s} = 
(1 - 1/\sqrt{\la})^{\la - \sqrt{\la}} \le e^{-\sqrt{\la} + 1}$, respectively.
As $\la\to\infty$, the largest of these values is $1/\la$, so we have 
$O\((l/\la)^l\) = O(1/\la)$ for $0\le l\le l_0$.
Thus the first sum is
$$\sum_{0\le k\le l} \left({l\over\la}\right)^k   = {1 + O(1/\la)\over 1 - l/\la}.$$
For the second sum we have
$$\sum_{0\le k\le l}{k\choose 2} \left({l\over\la}\right)^{k-1}   = {1 + O\(l^2(l/\la)^l\)\over (1 - l/\la)^3}.$$
The logarithm of $l^2(l/\la)^l$ has a non-negative second derivative for $l\ge 3$, so an argument similar to that used for the first sum shows that $O\(l^2(l/\la)^l\) = O(1/\la)$ for $0\le l\le l_0$.
Thus we have
$$\sum_{0\le k\le l}{k\choose 2} \left({l\over\la}\right)^{k-1}   = {1 + O(1/\la)\over (1 - l/\la)^3}$$
and the lower bound
$$D_l \ge {1 + O(1/\la)\over 1 - l/\la} - {1 + O(1/\la)\over \la(1 - l/\la)^3}. \eqno(3.3)$$

For an upper bound, we have
$$\eqalign{
l(l-1)\cdots(l-k+1) 
&\le l^k - \left(\sum_{0\le j\le k-1} j\right) l^{k-1} 
+ \left(\sum_{0\le i<j\le k-1} ij\right) l^{k-2} \cr
&\le l^k - {k\choose 2}l^{k-1} + {1\over 2}{k\choose 2}^2 l^{k-2} \cr
}$$
(because $\sum_{0\le i<j\le k-1} ij = \left(\left(\sum_{0\le j\le k-1} j\right)^2 - \sum_{0\le j\le k-1} j^2\right)\bigg/2
\le \left(\sum_{0\le j\le k-1} j\right)^2\big/2 = {k\choose 2}^2/2$).
Thus we have
$$D_l \le \sum_{0\le k\le l} \left({l\over\la}\right)^k  
- {1\over\la}\sum_{0\le k\le l} {k\choose 2}\left({l\over\la}\right)^{k-1}
+ {1\over2\la^2}\sum_{0\le k\le l} {k\choose 2}^2 \left({l\over\la}\right)^{k-2}.$$
For the third sum we have
$$\eqalign{
\sum_{0\le k\le l} {k\choose 2}^2 \left({l\over\la}\right)^{k-2}
&\le \sum_{k\ge 0} {k\choose 2}^2 \left({l\over\la}\right)^{k-2} \cr
&= O\left({1\over (1-l/\la)^5}\right). \cr
}$$
and thus the upper bound
$$D_l \le {1 + O(1/\la)\over 1 - l/\la} - {1 + O(1/\la)\over \la(1 - l/\la)^3}
 + O\left({1\over \la^2(1-l/\la)^5}\right).$$
 Combining this upper bound with the lower bound (3.3) yields
 $$D_l = {1 + O(1/\la)\over 1 - l/\la} - {1 + O(1/\la)\over \la(1 - l/\la)^3}
 + O\left({1\over \la^2(1-l/\la)^5}\right).$$
 
 To obtain $\Pr[L>l]$, we take the reciprocal of $D_l$:
 $$\eqalign{
 \Pr[L>l]
 &= \left({1 + O(1/\la)\over 1 - l/\la} - {1 + O(1/\la)\over \la(1 - l/\la)^3}
 + O\left({1\over \la^2(1-l/\la)^5}\right)\right)^{-1} \cr
  &= \(1 + O(1/\la)\)\;(1 - l/\la) \, \left(1 - {1 \over \la(1 - l/\la)^2}
 + O\left({1\over \la^2(1-l/\la)^4}\right)\right)^{-1} \cr
  &= \(1 + O(1/\la)\)\;(1 - l/\la) \, \left(1 + {1 \over \la(1 - l/\la)^2}
 + O\left({1\over \la^2(1-l/\la)^4}\right)\right) \cr
 &= \(1 + O(1/\la)\)\; \left((1 - l/\la) +  {1 \over \la(1 - l/\la)}  
 + O\left({1\over \la^2(1-l/\la)^3}\right)\right). \cr
 }$$
Observing that $O(1/\la)\,(1-l/\la) = O(1/\la)$ and
$O(1/\la)/\la(1-l/\la) = O\(1/\la^2(1-l/\la)^3\)$, we obtain (3.2).

We turn now to the tail.
We shall establish the estimate
$$\Pr[L>l] = O(e^{-\la} \, \la^l / l!) \eqno(3.4)$$
for $l\ge \la - s$, where $s = \sqrt{\la}$.
To obtain an upper bound on $\Pr[L>l]$, we obtain a lower bound on $D_l$.
We have
$$\eqalignno{
D_l
&= \sum_{0\le k\le l} {l! \over (l-k)! \, \la^k} \cr
&\ge {l! \over \lfloor \la-s\rfloor! \, \la^{l-\lfloor \la-s\rfloor}} + \cdots
+ {l! \over \lfloor \la-2s\rfloor! \, \la^{l-\lfloor \la-2s\rfloor}}, &(3.5)\cr
}$$
because $l-\lfloor \la-s\rfloor \ge l- (\la-s) \ge 0$ by assumption and 
$\lfloor \la-2s\rfloor \ge 0$ for all sufficiently large $\la$.
There are $\lfloor \la-2s\rfloor - \lfloor \la-2s\rfloor + 1 \ge s$ terms
in the sum (3.5).
Furthermore, the smallest of these terms is the last, because its denominator contains factors of 
$\la$ where the preceding terms contain factors smaller than $\la$.
Thus we have
$$D_l \ge {s \, l! \over \lfloor \la-2s\rfloor! \, \la^{l-\lfloor \la-2s\rfloor}}.$$
For the factorial in the denominator of this bound, we shall use the estimate
$n! \le e \, \sqrt{n} \, e^{-n} \, n^n$, which holds for all $n\ge 1$
(because the trapezoidal rule underestimates the integral $\int_1^n \log x\,dx$ of the concave function $\log x$).
This estimate yields
$$D_l \ge {s \, l! \, e^{\lfloor \la-2s\rfloor} \over 
e\, \sqrt{\lfloor \la-2s\rfloor} \, \lfloor \la-2s\rfloor^{\lfloor \la-2s\rfloor} 
\, \la^{l-\lfloor \la-2s\rfloor}}. \eqno(3.6)$$
We have
$$e^{\lfloor \la-2s\rfloor} \ge e^{\la-2s-1},$$
$$\eqalign{
\lfloor \la-2s\rfloor^{\lfloor \la-2s\rfloor} 
&\le (\la-2s)^{\lfloor \la-2s\rfloor} \cr
&= \la^{\lfloor \la-2s\rfloor} \, (1-2s/\la)^{\lfloor \la-2s\rfloor} \cr
&\le \la^{\lfloor \la-2s\rfloor} \, (1-2s/\la)^{\la-2s-1} \cr
&\le \la^{\lfloor \la-2s\rfloor} \, e^{(-2s/\la)(\la-2s-1)} \cr
&\le \la^{\lfloor \la-2s\rfloor} \, e^{-2s+4s^2/\la + 1} \cr
&\le \la^{\lfloor \la-2s\rfloor} \, e^{-2s+5} \cr
}$$
and
$$\sqrt{\lfloor \la-2s\rfloor} \le s.$$
Substituting these bounds into (3.6) yields
$$D_l \ge {l! \, e^\la \over e^7 \, \la^l}.$$
Taking the reciprocal of this bound yields (3.4).

We shall now use (3.2) and (3.4) to prove (1.8).
We write
$$\eqalign{
\De_m(l) 
&= l^m - (l-1)^m \cr
&= m\,l^{m-1} + O(l^{m-2}) \cr
}$$
for the backward differences of the $m$-th powers of $l$.
Then partial summation yields
$$\eqalignno{
\Ex[L^m]
&= \sum_{l\ge 0} l^m \, \Pr[L=l] \cr
&= \sum_{l\ge 0} \De_m(l) \, \Pr[L>l] \cr
&= \sum_{l\ge 0} m\, l^{m-1} \, \Pr[L>l] + O\left(\sum_{l\ge 0} l^{m-2} \, \Pr[L>l]\right) &(3.7)\cr
}$$
This formula shows that we should evaluate sums of the form
$$U_n = \sum_{l\ge 0} l^n \, \Pr[L>l]. \eqno(3.7)$$
We shall show that
$$U_n = {\la^{n+1} \over (n+1)(n+2)} + {\la^n \, \log \la \over 2} + O(\la^n). \eqno(3.8)$$
Substitution of this formula into (3.7) will then yield (1.8).

We shall break the range of summation in (3.7) at $l_0 = \la - s$, where $s=\sqrt{\la}$,
using (3.2) for $0\le l\le l_0$ and (3.4) for $l>l_0$.
Summing the first term in (3.2),
we have
$$\eqalignno{
\sum_{0\le l\le l_0} l^n (1 - l/\la)
&= {1\over\la} \sum_{0\le l\le l_0} (\la\,l^n - l^{n+1}) \cr
&= {1\over\la} \left(\left({\la \, l_0^{n+1}\over n+1} + O(l_0^n)\right)
- \left({\la^{n+2} \over n+2} + O(l_0^{n+1})\right)\right) \cr
&= {1\over\la} \left(\left({\la \, (\la^{n+1} - (n+1)\la^n s)\over n+1} + O(\la^n)\right)
- \left({\la^{n+2}  - (n+2)\la^{n+1} s\over n+2} + O(\la^{n+1})\right)\right) \cr
&= {\la^{n+1} \over (n+1)(n+2)} + O(\la^n).  \cr
}$$
Summing the second term in (3.2), we have
$$\eqalignno{
\sum_{0\le l\le l_0} {l^n \over \la-l}
&= \sum_{s\le k\le\la} {(\la - k)^n \over k} \cr
&= \sum_{s\le k\le\la} \left({\la^n \over k} + O(\la^{n-1})\right) \cr
&= \la^n \log {\la\over s} + O(\la^n) \cr
&= {\la^n \log\la \over 2} + O(\la^n), \cr
}$$
where we have used $\sum_{1\le k\le n} 1/k = \log n + O(1)$.
Summing the third term in (3.2) of course yields $O(\la^n)$.
Summing the last term in (3.2), we have
$$\eqalignno{
\la \sum_{0\le l\le l_0} {l^n \over (\la - l)^3}
& = \la \sum_{s\le k\le\la} {(\la-k)^n \over k^3} \cr
& \le \la^{n+1} \sum_{s\le k\le\la} {1 \over k^3} \cr
& \le \la^{n+1} \sum_{k\ge s} {1 \over k^3} \cr
& = \la^{n+1} \left({2\over s^2} + O\left({1\over s^3}\right)\right) \cr
&= O(\la^n), \cr
}$$
where we have used $\sum_{k\ge n} 1/k^3 = 2/n^2 + O(1/n^3)$.
Combining these estimates, we obtain
$$\sum_{0\le l\le l_0} l^n \, \Pr[L>l] 
= {\la^{n+1} \over (n+1)(n+2)} + {\la^n \, \log \la \over 2} + O(\la^n). \eqno(3.9)$$
Finally, summing (3.4) we have
$$\eqalign{
\sum_{l>l_0} {l^n \, e^{-\la} \, \la^l \over l!} 
&\le \sum_{l\ge 0} {l^n \, e^{-\la} \, \la^l \over l!}  \cr
&= O(\la^n), \cr
}$$
because the summation on the right-hand side is the $n$-th moment of a Poisson random variable with mean $\la$, which is a polynomial of degree $n$ in $\la$.
Thus
$$\sum_{l > l_0} l^n \, \Pr[L>l] = O(\la^n).$$
Combining this estimate with (3.9) yields (3.8) and completes the proof of (1.8).

Our final goal is to study the distribution of the random variable $I$, defined as the index of the first vacant waiting station $\calW_I$ found by a customer in the $M/M/1$ system, when the server, upon becoming free when at least one customer is waiting, serves the customer at the first occupied
waiting station.
We shall show that the distribution of $I$ is given by
$$\Pr[I > i] = {(1-\la) \, \la^{i+1} \over 1 - \la^{i+1}}. \eqno(3.10)$$
The event $I > i$ is simply the event that the first $i$ waiting stations $\calW_1,\ldots,\calW_i$ are all occupied; we have $I=0$ when the server is idle.
Letting the random variable $J$ denote the number of customers in the system, as before,
we observe that the event $J > i$ is necessary for the event $I > i$:
if stations $\calW_1,\ldots,\calW_i$ are occupied and the server is busy, there are at least $i+1$
customers in the system.
Thus we can write 
$$\eqalign{
\Pr[I > i] 
&= \sum_{j>i} \Pr[I > i \mid J=j] \, \Pr[J=j]. \cr
}$$ 
We shall show that
$$\Pr[I > i \mid J=j] =  {1-\la  \over 1 - \la^{i+1}} \eqno(3.11)$$
for all $j>i$.
Since $\Pr[J > i] = \la^{i+1}$, it will then follow that
$$\eqalign{
\Pr[I > i]
&= {1-\la  \over 1 - \la^{i+1}} \sum_{j>i} \Pr[J=j] \cr
&= {1-\la  \over 1 - \la^{i+1}}  \Pr[J > i] \cr
&= {(1-\la) \, \la^{i+1} \over 1 - \la^{i+1}},   \cr
}$$
confirming (3.10)

Before proving (3.11) in the general case, it will be helpful to consider two special cases.
If $i=0$, then the event $J > i$ is sufficient as well as necessary for the event $I > i$:
if the server is busy, an arriving customer must wait.
Thus $\Pr[I > i \mid J=j] = 1$ for all $j>0$, confirming (3.11) in this case. 
For $i=1$, we assume that $J =j > 1$ and ask for the conditional probability that $\calW_1$ is occupied.
If the current arrival occurs at time $t_0$, we consider the latest transition in the embedded Markov chain for $J$ that precedes $t_0$.
Suppose this previous transition occurs at time $t_1$.
If this previous transition was an arrival, then $\calW_1$ will be occupied by it (if it was not already occupied) and thus will be occupied at $t_0$.
If on the other hand this previous transition was a departure, then $\calW_1$ will be vacated by it (if it was not already vacant) and thus will be vacant at $t_0$.
Thus we must determine the probability that this previous transition was an arrival.
We claim that the previous transition was an arrival is $q = 1/(1+\la)$ and the probability that it was a departure is $p = \la/(1+\la)$.
To prove this claim, we note that the Markov chain for $J$ is {\it reversible}; that is,
if a movie is made of its transitions, the movie run backward is stochastically indistinguishable from the movie run forward.
(Reversibility follows from the fact that in this Markov chain, transitions occur only between  adjacent states; that is, $J$ is incremented by an arrival and decremented by a departure.)
The previous transition was an arrival if and only if it appears as a departure when the movie is run backward, and by reversibility this event occurs with probability $q = 1/(1+\la)$ provide that the transition does not involve the state $J=0$.
(When $J=0$, the next transition is an arrival with probability one, rather than with probability 
$p = \la/(1+\la)$.)
But since $J>1$ at $t_0$, we have $J>0$ at $t_1$.
This proves our claim.
We therefore have $\Pr[I > i \mid J=j] = q = 1/(1+\la) = (1-\la)/(1-\la^2)$, again confirming (3.11).

We are now ready to prove (3.11) in the general case.
We assume that $J =j> i$ and ask for the conditional probability that 
$\calW_1, \ldots, \calW_i$ are all occupied.
To determine whether this event occurs, we shall again trace backward in time through the transitions preceding the current arrival at time $t_0$.
In this case we may have to trace back through arbitrarily many transitions.
As we trace backward, we keep track of the difference $d(t)$ between the number of 
arrivals and the number of departures in the interval $[t,t_0)$.
We shall continue tracing as long as $-1<d(t)<i$, stopping at the latest time $t_1$ such that
$d(t_1) = -1$ or $d(t_1)=i$.

First, we claim that if $d(t_1)=i$, then $\calW_1, \ldots, \calW_i$ are all occupied at time $t_0$.
To prove this claim, we match each departure in $[t_1,t_0)$ with a later arrival ``like parentheses''.
That is, we associate each departure in this interval with a left parenthesis and each arrival with a right parenthesis.
Since $d(t)\ge 0$ for $t_1 \le t \le t_0$, there are at least as many right parentheses as left parentheses in any suffix of the resulting string, so all the left parentheses can be matched to right parentheses, leaving $d(t_1)=i$ right parentheses unmatched.
We thus have $i$ unmatched arrivals.
These arrivals occupy the stations $\calW_1, \ldots, \calW_i$ (if they were not already occupied),
and any stations that are vacated by subsequent departures are reoccupied by the matching arrivals, so these stations all remain occupied at time $t_0$.
This proves our first claim.

Next, we claim that if $d(t_1)=-1$, then at least one of the stations $\calW_1, \ldots, \calW_i$ is vacant at time $t_0$.
To prove this claim, we  define $m = \max\{d(t) ; t_1\le t\le t_0\}$.
We have $m<i$.
We define $t_2 = \max\{t: t_1<t\le t_0 \hbox{\ and\ } d(t)=m\}$.
Since $d(t)\le m$ for $t_1\le t\le t_2$, we can match arrivals with later departures in this interval,
leaving at least $m+1$ departures unmatched.
These departures will vacate stations $\calW_1, \ldots, \calW_{m+1}$ (if they were not already vacant), and any of these stations that are occupied by subsequent arrivals in this interval will be revacated by the matching departures, so that these stations will all be vacant at time $t_2$.
For each of these $m+1$ stations $\calW_j$, let $e(j)$ denote the difference between the number of arrivals that occupy $\calW_j$ and the number of departures that vacate $\calW_j$.
We have $0\le e(j)\le 1$, because arrivals that occupy a given station alternate with departures that vacate it.
Since $d(t_2) = m$, we have $\sum_{1\le j\le m+1} e(j) = m$.
It follows that $e(j)=0$ for at least one value of $j$, and so $\calW_j$ remains vacant at time $t_0$
for this value of $j$.
This proves our second claim.

Finally, we claim that as we trace backward, the probability that the next transition is an arrival is always $q = 1/(1+\la)$, and the probability that the next transition is a departure is always 
$p=\la/(1+\la)$.
To prove this claim, we observe that there are at least $i+1$ customers in the system at time $t_0$, and thus at least $i+1-d(t)$ customers in the system at any time $t$ for $t_1< t <t_0$.
Furthermore $i+1-d(t)\ge 2$, because $d(t)<i$ for
$t_1 < t < t_0$.
Thus there are at least two customers in the system at every time $t\in(t_1,t_0)$.
It follows that as we trace backward, none of the transitions we encounter involve the state $J=0$.
Thus as we trace backward, the probability that the next transition is an arrival is always
$q = 1/(1+\la)$, and the probability that the next transition is a departure is always $p=\la/(1+\la)$.
This proves our third claim.

These three claims show that the process of determining whether $I > i$ (given that $J =j> i$)
is as shown in Figure 3.
\medskip
$$\xymatrix{
&&\ar[d]^-{\text{start}} \\
 & *+++[o][F=]{\scriptstyle I\le i}  & \state{\scriptstyle 0} \ar@/^/[r]^-q \ar[l]^-p&\cdots \ar@/^/[r]^-q \ar@/^/[l]^-p & \state{\scriptstyle d} \ar@/^/[r]^-q \ar@/^/[l]^-p & \cdots \ar@/^/[r]^-q \ar@/^/[l]^-p &\state{\scriptstyle i-1}\ar[r]^-q  \ar@/^/[l]^-p&*+++[o][F=]{\scriptstyle I > i}
}$$
\medskip
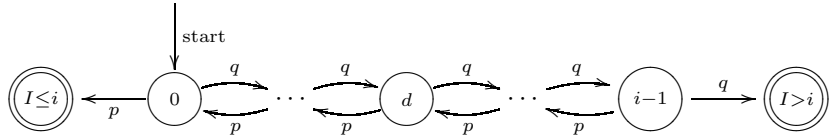
\centerline{Figure 3. State transition diagram for determining whether $I\le i$ or $I > i$
(given that $J = j > i$).}
\medskip
Comparison with Figure 2 shows that this process is the same at the process of determining whether $K > k$, except that the roles of $p$ and $q$ are exchanged.
This exchange is equivalent to the substitution of $1/\la$ for $\la$, so $\Pr[I > i \mid J=j]$ is obtained by making this substitution in the expression (2.2) for $\Pr[K > k]$;
$$\eqalign{
\Pr[I > i \mid J=j]
&= {(1 - 1/\la) \,(1/\la)^i \over \(1 - (1/\la)^{i+1}\)} \cr
&= {1 - \la \over 1 - \la^{i+1}}. \cr
}$$
This again confirms (3.11), which completes the proof of (3.10).
It follows that the moments of $I$ differ from those of $K$ by a factor of $\la = 1 - (1-\la)$.
This fact allows asymptotic expansions for the moments of $I$ to be obtained from those of $K$,
with the result that the leading terms are the same.
\sk

\heading{4. Conclusion}

We observe that the search algorithm described in the preceding section  for the $M/M/1$ system  defines a deterministic service discipline distinct from both first-come-first-served and last-come-first-served.
It would be of interest to determine the distribution of the waiting time $W$ experienced by a customer for this discipline, or even the moments of this distribution.
We hope to address this question in a future paper.
\sk

\heading{5. Acknowledgment}

The research reported here was supported
by Grant CCF  0917026 from the National Science Foundation.
\sk

\heading{6. References}

\ref C1; E. G. Coffman, Jr., T. T. Kadota and L. A. Shepp;
``A Stochastic Model of Fragmentation in Dynamic Storage Allocation'';
SIAM J. Comput.; 14:2 (1985) 416--425.

\refbook C2;  J. W. Cohen;
The Single Server Queue (revised edition);
North-Holland, Amsterdam, 1982.

\refbook C3; L. Comtet;
Advanced Combinatorics:
The Art of Finite and Infinite Expansion;
D.~Reidel Publishing Co., Dortrecht, 1974. 

\ref E1; S. Egger (n\'{e} Endres) and F. Steiner;
``A New Proof of the Vorono\"{\i} Summation Formula'';
J. Phys.\ A: Math.\ Theor.; 44 (2011) 225302 (11 pp.).

\ref E2; S. Endres and F. Steiner;
``A Simple Infinite Quantum Graph'';
Ulmer Seminare Funktionalanalysis und Differentialgleichungen'';
14 (2009) 187--200.

\refbook F; W. Feller;
An Introduction to Probability Theory and Its Applications (3rd edition);
John Wiley \& Sons, New York, 1968.

\refbook G1; G. Gasper and M. Rahman;
Basic Hypergeometric Series;
Cambridge University Press, Cambridge, 1990.

\ref G2; E. Godehardt and J. Jaworski;
``On the Connectivity of a Random Interval Graph'';
Random Struct.\ Alg.; 9 (1996) 137--161.

\refbook H; G. H. Hardy and E. M. Wright;
Introduction to the Theory of Numbers (5th edition);
Clarendon Press, Oxford, 1979.

\refbook N; G. F. Newell;
The $M/M/\infty$ Service System with Ranked Servers in Heavy Traffic;
Springer-Verlag, Berlin, 1984.

\ref P; N. Pippenger;
``Random Interval Graphs'';
Random Struct.\ Alg.; 12 (1998) 361--380.

\refbook R1; J. Riordan;
Stochastic Service Systems;
John Wiley and Sons, New York, 1962.

\refbook R2; S. Roman;
The Umbral Calculus;
Academic Press, New York, 1984.

\ref S1; E. R. Scheinerman; 
``Random Interval Graphs'';
Combinatorica; 8 (1988) 357--371.

\ref S2; E. R. Scheinerman; 
``An Evolution of Interval Graphs'';
Discrete Math.; 82 (1990) 287--302.

\refbook W; E. T. Whittaker and G. N. Watson;
A Course of Modern Analysis (4th edition);
Cambridge University Press, London, 1927.

\refinbook Z; D. Zagier;
``The Mellin Transformation and Other Useful Analytic Techniques'';
in: E.~Zeidler; Quantum Field Theory: I, Basics in Mathematics and Physics;
Springer-Verlag, Berlin, 2006, pp.~305--323.

\bye